\documentclass{amsart}
\font\fiverm=cmr5

\let\<=\langle
\let\>=\rangle
\let\sse=\subseteq
\let\noi=\noindent
\let\veps=\varepsilon
\let\limply=\Longrightarrow

\def\0{\{0\}}
\def\span{{\kern.5pt{\rm span}\kern1pt}}
\def\smallfrac#1#2{{\textstyle{\frac{#1}{#2}}}}
\def\conv{{\;\longrightarrow\;}}
\def\wconv{{{\buildrel_{\scriptstyle w}\over\conv}}}
\def\sconv{{{\buildrel_{\scriptstyle s}\over\conv}}}
\def\uconv{{{\buildrel_{\scriptstyle u}\over\conv}}}
\def\sslash{\hbox{{\fiverm/}}}
\def\notconv{{{\conv\kern-13pt\slash}\kern9pt}}
\def\notuconv{{{\uconv\kern-13pt\sslash}\kern9pt}}
\def\notsconv{{{\sconv\kern-13pt\sslash}\kern9pt}}
\def\notwconv{{{\wconv\kern-13pt\sslash}\kern9pt}}

\def\B{{\mathcal B}}

\def\Oe{{\mathcal O}}
\def\R{{\mathcal R}}
\def\X{{\mathcal X}}
\def\BX{{\B[\X]}}

\def\CC{{\mathbb C\kern.5pt}}
\def\DD{{\mathbb D\kern.5pt}}
\def\FF{{\mathbb F\kern.5pt}}
\def\JJ{{\mathbb J\kern.5pt}}
\def\RR{{\mathbb R\kern.5pt}}
\def\TT{{\mathbb T\kern.5pt}}

\def\jj{{\scriptstyle\JJ}}

\newsymbol\varnothing 203F
\let\void=\varnothing

\def\matrix#1{\null\,\vcenter{
              \normalbaselines\mathsurround=0pt\ialign{
              \hfil $##$
              \hfil && \quad
              \hfil $##$
              \hfil \crcr
              \mathstrut \crcr
              \noalign{\kern-\baselineskip}#1 \crcr
              \mathstrut \crcr
              \noalign{\kern-\baselineskip} \crcr }}\,}

\def\newmatrix#1{\null\,\vcenter{
		\baselineskip=8pt\mathsurround=-0pt\ialign{
		\hfil ${##}$
		\hfil &&
		\hfil ${##}$
		\hfil \crcr
		\mathstrut \crcr
		\noalign{\kern-\baselineskip}#1 \crcr
		\mathstrut \crcr
		\noalign{\kern-\baselineskip} \crcr }}\!}

\def\smallmatrix#1{\null\,\vcenter{
		\baselineskip=8pt\mathsurround=-0pt\ialign{
		\hfil ${\scriptstyle##}$
		\hfil &&
		\hfil ${\scriptstyle##}$
		\hfil \crcr
		\mathstrut \crcr
		\noalign{\kern-\baselineskip}#1 \crcr
		\mathstrut \crcr
		\noalign{\kern-\baselineskip} \crcr }}\!}

\begin{document}

\vglue-35pt
\hfill{\it Rendiconti del Circolo Matematico di Palermo}\/,
{\bf 73}(8) (2024) 3217--3228

\vglue30pt
\title{Weak Stability and Quasistability}
\author{C.S. Kubrusly}
\address{Catholic University of Rio de Janeiro, Brazil}
\email{carlos@ele.puc-rio.br}
\author{P.C.M. Vieira}
\address{National Laboratory for Scientific Computation,
         Petr\'opolis, RJ, Brazil}
\email{paulocm@lncc.br}
\renewcommand{\keywordsname}{Keywords}
    \date{February 12, 2024}

\begin{abstract}
It is known that weak l-sequential supercyclicity implies weak quasistability,
and it is still unknown whether weak l-sequential supercyclicity implies
weak stability, much less whether weak supercyclicity implies weak stability
(although it has been known for a long time that strong supercyclicity implies
strong stability)$.$ It is shown that weak l-sequential supercyclicity implies
weak stability under the assumption of boundedly spaced subsequences$.$ It is
given an example of a power bounded weakly quasistable operator that is not
weakly stable.
\end{abstract}

\maketitle

\vskip-15pt\noi
\section{\bf Introduction}

Let $\X$ be a normed space, and let $\X^*$ be its dual$.$ An $\X$-valued
sequence $\{x_n\}$ converges weakly if the scalar-valued sequence $\{f(x_n)\}$
converges for every ${f\in\X^*}\!.$ Take an operator $T$ on $\X.$ It is weakly
stable if $\lim_n|f(T^nx)|\!=0$, and weakly quasistable if
$\liminf_n\!|f(T^nx)|=0$, for every ${x\in\X}$, for every ${f\in\X^*}\!.$ It
is clear that weak stability implies weak quasistability$.$ It is shown that
the converse fails by exhibiting a power bounded weakly quasistable operator
that is not weakly stable (\hbox{Proposition 4.3}).

\vskip6pt
Let $A$ be any nonempty set$.$ A subsequence $\{a_{n_k}\}$ of an $A$-valued
sequence $\{a_n\}$ is boundedly spaced if ${\sup_k(n_{k+1}-n_k)<\infty}.$ A
central result in this paper shows that if a weakly quasistable operator has
a boundedly spaced subsequence of weak quasistability for every vector, then
it is weakly stable (Theorem 5.3).

\vskip6pt
Let ${\Oe_T([y])}$ be the projective orbit of a nonzero vector ${y\in\X}$
under an operator $T.$ The operator $T$ is weakly l-sequentially supercyclic
if every ${x\in\X}$ is the weak limit of an $\Oe_T([y])$-valued sequence$.$
It is unknown whether weak l-sequential supercyclicity implies weak
stability$.$ An application of Theorem 5.3 gives a sufficient condition for
weak l-sequential supercyclicity to imply weak stability
\hbox{(Theorem 6.2)}$.$

\vskip6pt
On the other hand, it is also given a necessary condition for weak stability
of weakly l-sequentially supercyclic operators (Theorem 6.3), which together
with Theorem 6.2, shows what happens when weak l-sequential supercyclicity
holds without boundedly spaced subsequences (Corollary 6.4).

\vskip6pt
The paper is organised as follows$.$ Basic notation and terminology are
\hbox{summarised} in Section 2$.$ The notions of uniform, strong, and weak
stabilities, and their quasi\-stability counterparts, are considered in
Section 3$.$ It is shown in Section 4 that the notions of stability and
quasistability coincide for the strong and uniform cases, but not for the weak
case$.$ Boundedly spaced subsequences are considered in \hbox{Section} 5, and
an application to weak l-sequential supercyclicity is considered in
\hbox{Section 6}.

\section{\bf Basic notation and terminology }

Throughout this paper, all linear spaces are over the same scalar field $\FF$,
which is either $\RR$ or $\CC$, and $\X$ denotes an infinite-dimensional
normed space$.$ Particular cases of inner product, Banach, and Hilbert spaces
are discussed accordingly$.$ Let ${\BX=\B[\X,\X]}$ stand for the collection of
all bounded linear transformations of $\X$ into itself$.$ That is, $\BX$
stands for the normed algebra of all operators on $\X.$ The linear manifold
${\R(T)=T(\X)}$ of $\X$ is the range of $T$ in $\BX.$ The Banach~space
${\X^*\!=\B[\X,\FF]}$ is the dual of $\X$, and $T^*\kern-1pt$ in $\B[\X^*]$
denotes the normed-space adjoint of $T$ in $\BX$ (i.e., ${T^*\!f=f\kern1ptT}$
for every ${f\in\X^*}).\!$ We use the same notation for the Hilbert-space
adjoint ${T^*\kern-1pt\in\BX}$ of ${T\kern-1pt\in\BX}$ on a Hilbert space
$\X$, where the concepts of dual and adjoint are shaped after the Riesz
Representation Theorem in Hilbert space --- in a Hilbert space setting,
adjoints are always Hilbert-spaces's.) An $\X$-valued sequence $\{x_n\}$
converges weakly to ${x\in\X}$ if ${\lim_nf(x_n)\to f(x)}$ for every
${f\in\X^*}$, equivalently, ${\lim_nf(x_n-x)\to0}$ for every ${f\in\X^*}.$ 
(This definition is standard --- see, e.g., \cite[Definition 1.13.2]{Meg})$.$
Alternative and usual notation for weak convergence that will be used here:
${x_n\wconv x}$ and $x={w\hbox{\,-}\lim_n x_n}.$ If $\X$ is a Hilbert space
with inner product ${\<\,\cdot\,;\cdot\kern1pt\>}$, then (again, by the Riesz
Representation Theorem in Hilbert space) weak convergence of $\{x_n\}$ to $x$
means ${\<x_n\!-x\,;z\>\to0}$ for \hbox{every} ${z\in\X}.$ Norms on $\X$ and
the induced uniform norm on $\B[\X]$ (and the norm on $\X^*$) will be denoted
by the same symbol ${\|\,\cdot\,\|}.$ An $\X$-valued sequence $\{x_n\}$ is
strongly convergent (or converges in the norm topology) to ${x\in\X}$ if
${\|x_n\!-x\|\to0}.$ Alternative notation: ${x_n\!\sconv x}$; it is clear
that strong convergence implies weak convergence$.$ For an arbitrary operator
${T\kern-1pt\in\BX}$, take its power sequence $\{T^n\}.$ According to the
above definitions of convergence in $\X$, an operator $T$ is said to be weakly
or strongly stable if ${T^nx\wconv0}$ or ${T^nx\sconv0}$ for every ${x\in\X}$,
respectively$.$ An operator $T$ is power bounded if
${\sup_n\kern-1pt\|T^n\|<\infty}$ (equivalently, if
${\sup_n\kern-1pt\|T^nx\|<\infty}$ for every ${x\in\X}$ by the
Banach--Steinhaus Theorem, if \hbox{$\X$ is a Banach space).}

\section{\bf Stability and quasistability}

Take ${T\kern-1pt\in\BX}$ on a normed space $\X$ and consider the $\BX$-valued
power sequence $\{T^n\}.$ The operator $T$ is uniformly stable if $\{T^n\}$
converges uniformly (i.e., if the sequence $\{T^n\}$ converges in the norm
topology on $\BX$) to the null operator$.$ Thus an operator $T$ is uniformly,
strongly, or weakly stable if
$$
{\lim}_n\|T^n\|=0,
\quad
\hbox{or}
\quad
{\lim}_n\|T^nx\|\to0
\;\;\hbox{for every}\;\;
x\in\X,
\quad
\hbox{or}
$$
\vskip-2pt\noi
$$
{\lim}_n|f(T^nx)|=0
\;\;\hbox{for every}\;\;
x\in\X,
\;\;\hbox{for every}\;\;
f\in\X^*,
$$
\vskip2pt\noi
respectively$.$ Alternative notations for strong and weak stability were
given in the previous section and, accordingly, if $\X$ is a Hilbert space,
then weak stability for $T$ means ${\<T^nx\,;z\>\to0}$ for every ${x,z\in\X}$
--- if the Hilbert space $\X$ is complex, this is equivalent to
${\<T^nx\,;x\>\to0}$ for every ${x\in\X}.$ It is clear that uniform stability
implies strong stability, which implies weak stability, which in turn implies
power boundedness (if $\X$ is a Banach space), and the converses fail$.$ We
say that an operator ${T\kern-1pt\in\BX}$ on a normed space $\X$ is
{\it uniformly quasistable}\/, or {\it strongly quasistable}\/, or {\it weakly
quasistable}\/ if the above limits are replaced by limits inferior$.$ That is,
if
$$
{\liminf}_n\|T^n\|=0,
\quad
\hbox{or}
\quad
{\liminf}_n\|T^nx\|\to0
\;\;\hbox{for every}\;\;
x\in\X,
\quad
\hbox{or}
$$
\vskip-2pt\noi
$$
{\liminf}_n|f(T^nx)|=0
\;\;\hbox{for every}\;\;
x\in\X,
\;\;\hbox{for every}\;\;
f\in\X^*,
$$
\vskip2pt\noi
respectively$.$ It is clear that each form of stability implies its respective
form of quasistability, and each form of quasistability implies the existence
of a subsequence of the power sequence that converges to zero$.$ In other
words,
\vskip6pt\noi
(i)
uniform quasistability means that there exists a subsequence $\{T^{n_k}\}$ of
$\{T^n\}$ that converges in $\BX$ to the null operator, that is,
${\|T^{n_k}\|\to0}$;
\vskip6pt\noi
(ii)
strong quasistability means that for every ${x\in\X}$ there exists a
subsequence $\{n_k\}=\{n_k(x)\}$ of the sequence of the positive integers
(which depends on $x$) such that the $\X$-valued sequence $\{T^{n_k}x\}$
converges in the norm topology to the origin of $\X.$  That is,
${T^{n_k}x\sconv0}$ for every ${x\in\X}$, for some subsequence $\{T^{n_k}\}$
of $\{T^n\}$ --- equivalent notation: ${\|T^{n_k}x\|\to0}$ for every
${x\in\X}$;
\vskip4pt\noi
(iii)
weak quasistability means that for each ${x\in\X}$ there is a subsequence
$\{n_k\}=\{n_k(x)\}$ of the sequence of the positive integers (which depends
on $x$) such that ${T^{n_k}x\wconv0}.$ Equivalently, for each vector
${x\in\X}$ there exists a subsequence $\{T^{n_k}\}=\{T^{n_k(x)}\}$ of
$\{T^n\}$ such that the scalar-valued sequences $\{f(T^{n_k}x)\}$ converge to
zero for every functional ${f\in\X^*}\!.$ We summarise this by saying that
there is a subsequence $\{T^{n_k}\}$ of $\{T^n\}$ such that
${|f(T^{n_k}x)|\to0}$ for every ${x\in\X}$, for \hbox{every ${f\in\X^*}\!$.}

\vskip6pt
We will be particularly concerned with weak quasistability$:$ for every $x$
in $\X$, ${T^{n_k}x\wconv0}$ for some subsequence $\{T^{n_k}\}$ of $\{T^n\}.$
Thus an operator ${T\kern-1pt\in\BX}$ is weakly quasistable if, for every
${x\in\X}$, there exists a subsequence $\{T^{n_k}\}$ of $\{T^n\}$ for which
${\lim_kf(T^{n_k}x)=0}$ (equivalently, ${\lim_k|f(T^{n_k}x)|=0}$) for every
${f\in\X^*}\!$ (where, for each ${x\in\X}$, the subsequence $\{T^{n_k(x)}\}$
depends on $x$ but not on $f).$

\section{\bf Auxiliary results}

To begin with, we verify that the notions of uniform and strong stabilities
coin\-cide with the motions of uniform and strong quasistabilities,
respectively$.$ The notion of weak stability, however, does not coincide
with the notion of weak quasi\-stability$.$ As we saw above, uniform and
strong stabilities imply uniform and strong quasistability, respectively and
trivially$.$ Thus we verify next only the converses for normed-space operators
(for a Hilbert-space version, see \cite[\hbox{proof of Theorem 3}]{KV1}).

\vskip6pt\noi
{\bf Proposition 4.1$.$}
{\it For every operator\/ $T$ on a normed space\/ $\X$, uniform quasistability
is equivalent to uniform stability$.$ Indeed,
$$
{\inf}_n\|T^n\|<1
\quad\;\limply\;\quad
{\lim}_n\|T^n\|=0,
$$
\vskip-4pt\noi
and hence}
$$
{\liminf}_n\|T^n\|=0
\quad\!\iff\quad\!
{\lim}_n\|T^n\|=0.
$$

\proof
Let $T$ be an operator on a normed space$.$ Recall that
${\lim_n\|T^n\|^{1/n}}$ exists as a real number and, by setting
${r(T)=\lim_n\|T^n\|^{1/n}}\!$, also recall that
$$
r(T)^n=r(T^n)\le\|T^n\|\le\|T\|^n
$$
for every nonnegative integer $n$ (see, e.g., \cite[Lemma 6.8]{EOT} --- in a
Banach space setting, which is not necessarily the case here, $r(T)$ coincides
with the spectral radius of $T\kern-1pt$, and uniform stability is equivalent
to ${r(T)<1}.$) However, the implication
$$
r(T)={\lim}_n\|T^n\|^{1/n}<1
\quad\limply\quad
{\lim}_n\|T^n\|\to0
$$
holds in a normed space as well (see, e.g.,
\cite[proof of Proposition 6.22]{EOT})$.$ Now, if ${\inf_n\|T^n\|<1}$, then
either there is an $n_0$ such that ${\|T^{n_0}\|<1}$, or
${\liminf_n\|T^n\|<1}$, which also ensures that there exists an $n_0$ such
that ${\|T^{n_0}\|<1}.$ Hence
$$
r(T)^{n_0}=r(T^{n_0})\le\|T^{n_0}\|<1,
$$
and therefore ${r(T)<1}$, which implies ${\lim_n\|T^n\|=0}.$ Thus
$$
{\inf}_n\|T^n\|<1
\quad\limply\quad
{\lim}_n\|T^n\|=0.
$$
So ${\liminf_n\|T^n\|=0}$ implies ${\lim_n\|T^n\|=0}$, leading to the stated
equivalence.                                                            \qed

\vskip6pt
Note that ${\lim_n\|T^n\|=0}$ naturally implies power boundedness for $T$.

\vskip6pt\noi
{\bf Proposition 4.2$.$}
{\it If\/ $T$ is a power bounded operator on a normed space\/ $\X$, then
strong quasistability is equivalent to strong stability}\/:
$$
{\liminf}_n\|T^nx\|=0
\;\;\hbox{\it for every}\;\;
x\in\X
\quad\iff\quad
{\lim}_n\|T^nx\|=0
\;\;\hbox{\it for every}\;\;
x\in\X.
$$

\proof
Take an arbitrary ${x\in\X}.$ If ${\liminf_n\|T^nx\|=0}$, then there exists
a subsequence $\{T^{n_k}x\}$ of $\{T^nx\}$ (which depends on ${x\in\X}$) for
which ${\lim_k\|T^{n_k}x\|=0}.$ Thus, if $T$ is power bounded,
$$
\|T^nx\|\le\|T^{n-n_k}\|\,\|T^{n_k}x\|\le{\sup}_n\|T^n\|\|T^{n_k}x\|
\quad\hbox{whenever}\quad
n\ge n_k.
$$
By the above inequality,
$$
{\lim}_k\|T^{n_k}x\|\!=\!0
\;\Leftrightarrow\;
{\limsup}_k\|T^{n_k}x\|\!=\!0
\;\Rightarrow\;
{\limsup}_n\|T^nx\|\!=\!0
\;\Leftrightarrow\;
{\lim}_n\|T^nx\|\!=0.
$$
Therefore, if $T$ is power bounded,
$$
{\liminf}_n\|T^nx\|=0
\;\;\hbox{\it for every $x\in\X$}
\;\;\Rightarrow
\;\;{\lim}_n\|T^nx\|=0
\;\;\hbox{\it for every $x\in\X$}.                              \eqno{\qed}
$$

\vskip6pt
However, although weak stability implies weak quasistability, the
converse fails$:$
$$
{\liminf}_n|f(T^nx)|=0
\;\;\;\hbox{for every}\;\;\;
x\in\X,
\;\;\;\hbox{for every}\;\;\;
f\in\X^*
$$
$$
{\kern4pt\not\kern-4pt\limply}\;\;
{\lim}_n|f(T^nx)|=0
\;\;\hbox{for every}\;\;
x\in\X,
\;\;\hbox{for every}\;\;
f\in\X^*\!.
$$
\vskip6pt

We show in Proposition 4.3 below that the Foguel operator, which is power
bounded and acts on a separable Hilbert space, is weakly quasistable but not
weakly stable$.$ Let $\X$ be a separable Hilbert space, and consider the
Foguel operator
$$
F=\big(\smallmatrix{S^* &          P      \cr
                    O   & S^{\phantom{*}} \cr}\big)
$$
acting on the separable Hilbert space ${\X\oplus\X}$ (where $\oplus$ means
orthogonal direct sum) \cite{Fog,Hal}$.$ Here $S$ is a unilateral shift of
multiplicity one on $\X$, which shifts an orthonormal basis $\{e_k\}_{k\ge0}$
for $\X$, and ${P\!:\X\!\to\X}$ is the orthogonal projection onto
$\R(P)={\span\{e_j\!:j\in\JJ\}^-}\!$ (closure of span of $\{e_j\!:j\in\JJ\}$),
where $\JJ$ is a sparse infinite set of positive integers with the following
property: if ${i,j\in\JJ}$ and ${i<j}$, then ${2i<j}$ (e.g.,
$\JJ={\{j\ge 1\!:\;j=3^k;\;k\ge 0\}}$, the set of all integral powers of 3).

\vskip6pt\noi
{\bf Proposition 4.3$.$}
{\it The Foguel operator is a power bounded operator that is weakly
quasistable but not weakly stable}\/.

\proof
Let $F$ be the Foguel operator defined above$.$ It is well known that $F$ is
power bounded \cite[p.791]{Hal} and not weakly stable (see, e.g.,
\cite[Remark 8.7]{MDOT})$.$ Note that
$$
F^n=\big(\smallmatrix{S^{*n} &          P_n     \cr
                      O      & S^{n\phantom{*}} \cr}\big)
\quad\hbox{for every}\quad
n\ge0,
$$
with ${P_n\:\X\!\to\X}$ given by ${P_{n+1}={\sum}^n_{i=0}S^{*\,n-i}PS^i}$ and
${P_0=O}$ (the null operator)$.$ Therefore, for every
${x=(x_1,x_2)\in\X\oplus\X}$ and every ${z=(z_1,z_2)\in\X\oplus\X}$,
$$
\<F^{n+1}x\,;z\>=\<S^{*n+1}x_1,\;z_1\>+\<P_{n+1}x_2\,;z_1\>+\<S^{n+1}x_2,z_2\>.
$$
Since $S$ is weakly stable (and so is $S^*$), we get
$$
{\liminf}_n|\<F^nx\,;z\>|={\liminf}_n|\<P_nx_2\,;z_1\>|            \eqno{(\S)}
$$
for every ${x=(x_1,x_2)\in\X\oplus\X}$, for every
${z=(z_1,z_2)\in\X\oplus\X}$, where
$$
\<P_{n+1}x_2\,;z_1\>={\sum}_k\<x_2;e_k\>\<P_{n+1}e_k\,;z_1\>,
$$
\kern-4pt\noi
with
$$
P_{n+1}e_k
={{\sum}^n_{i=0}S^{*\,n-i}Pe_{k+i}}
={{\sum}^{k+n}_{i=k}S^{*\,k+n-i}Pe_i}.
$$
\vskip4pt\noi
If ${i\not\in\JJ}$, then ${Pe_i=0}$; if ${i\in\JJ}$, then ${Pe_i=e_i}.$ Hence
$$
S^{*\,k+n-i}Pe_i
=\left\{\newmatrix{e_{2i-(k+n)}\quad & {\rm if}\:\;k+n\le2i,\;\;i\in\JJ, \cr
                                     &                                   \cr
                  {\kern-44pt0}      & {\rm{otherwise}}\,.               \cr}
 \right.
$$
Given a pair of nonnegative integers ${(k,n)}$, consider the set
$$
\jj_{_{\scriptstyle{k,n}}}
=\big\{j\in\JJ\:\;\;k\le j\le k+n\le 2j\big\}
=\JJ\cap\big[\max\{k,\smallfrac{k+n}{2}\}\,,\,k+n\big],
$$
which has at most one element (i.e., $\#{\jj_{k,n}\le1}).$ (Indeed, if
${i,j\in\jj_{k,n}}$ with ${i<j}$, then $2i<j$ (for ${i,j\in\JJ}$) and
${j\le k+n\le 2i}$, which is a contradiction.) Then
$$
P_{n+1}e_k
=\left\{\newmatrix{\sum_{i\in{\jj}_{k,n}}e_{2i-(k+n)}
                  \,=\,e_{2j-(k+n)}\quad & {\rm for}\;\;j\in\jj_{k,n}\;\;
                                     {\rm whenever}\;\;\jj_{k,n}\ne\void, \cr
                                         &                                \cr
                          {\kern-134pt0} & {\rm{otherwise}}\,.            \cr}
 \right.
$$
Now take an arbitrary pair ${(x_1,z_2)\in\X\oplus\X}.$ Since $P$ is an
orthogonal projection, $P_{n+1}$ is self-adjoint, thus we get
${\<P_{n+1}x_1\,;z_2\>}={\sum_k\<x_1\,;P_{n+1}e_k\>\<e_k\,;z_2\>}$, and so
$$
\<P_{n+1}x_1\,;z_2\>
=
\left\{\newmatrix{\sum_k\<x_1\,;e_{2j-(k+n)}\>\<e_k\,;z_2\>\quad
                                     & {\rm for}\;\;j\in\jj_{k,n}\;\;
                                     {\rm whenever}\;\;\jj_{k,n}\ne\void, \cr
                                     &                                \cr
                      {\kern-116pt0} & {\rm{otherwise}}\,.            \cr}
 \right.
$$
Then we can take an infinite subsequence $\{P_{n_j}\}$ of $\{P_n\}$ for which
${\<P_{n_j}x_1\,;z_2\>=0}$ for every $j$, and so
${\lim_j\<P_{n_j}x_1\,;z_2\>=0}$, for every ${x_1,z_2\in\X}\!.$ For instance,
with ${n_j=2j+1}$ for ${j\in\JJ}$ we get ${\liminf_n|\<P_nx_1\,;z_2\>|=0}$ for
every ${x_1,z_2\in\X}$, which implies ${\liminf_n|\<F^nx,z\>|=0}$ for every
${x,z\in\X\oplus\X}$ by $(\S)$.                                           \qed

\vskip6pt
Observe that $\{n_j\}$ is a rather sparse subsequence of the nonnegative
integers, whose increments ${n_{j+1}\!-n_j}$ increase unboundedly as
${j\to\infty}$ (by the construction of the set $\JJ).$ Such a property was
fundamental for proving that $F$ is not weakly stable (see, e.g.,
\cite[Proposition 8.3 and Remark 8.7]{MDOT}).

\section{\bf Boundedly spaced subsequences}

Let $\!\{n\}_{n\ge1}\!$ stand for the self-indexing of the set of all positive
integers equipped with the natural order and regard it as the
positive-integer-valued identity sequence$.$ A subsequence $\{n_k\}_{k\ge1}$
of $\{n\}_{n\ge1}$ is of {\it bounded increments}\/ (or has bounded gaps) if
${\sup_{k\ge1}(n_{k+1}-n_k)<\infty}.$ Let $A$ be an arbitrary nonempty set$.$
We say that a subsequence $\{a_{n_k}\}$ of an arbitrary $A$-valued sequence
$\{a_n\}$ is {\it boundedly spaced}\/ if it is indexed by a subsequence
$\{n_k\}_{k\ge1}$ of bounded increments (i.e., $\{a_{n_k}\}$ is boundedly
spaced if ${\sup_{k\ge1}{(n_{k+1}-n_k})<\infty}).$ These notions (bounded
increments and boundedly spaced) have been applied, in the present context,
in \cite{KV1,KV2}$.$ The next proposition will be required in the sequel$.$
(A standard result, for the particular case of scalar sequences, with
an easy proof, which we include for sake of \hbox{completeness}.)

\vskip6pt\noi
{\bf Proposition 5.1.}
{\it If\/ $\{\alpha_{n_k}\}$ is a boundedly spaced subsequence of an\/
$\FF$-valued se\-quence\/ $\{\alpha_n\}$, then\/
${\alpha_{n_k+j}\!\to \alpha}$ as\/ ${k\kern-1pt\to\kern-1pt\infty}$ for
every\/ ${j\kern-1pt\ge\kern-1pt1}$ implies\/ ${\alpha_n\!\to \alpha}$ as}\/
${n\kern-1pt\to\kern-1pt\infty}$.

\proof
For each ${j\ge1}$, suppose the sequence $\{\alpha_{n_k+j}\}$ converges to
$\alpha$ as ${k\to\infty}.$ In addition, suppose $j$ lies in a bounded set
$J$, and take an arbitrary ${\veps>0}.$ Thus there is an integer $k_\veps$
such that for all ${j\in[1,\max j]}$, with ${\max{j}=\max\{j:\! j\in J\}}$,
$$
k\ge k_\veps\quad\limply\quad|\alpha_{n_k+j}-\alpha|<\veps,
$$
(Indeed, for any ${\veps>0}$ and each $j$ take $k_{\veps,j}$ such that if
${k\ge k_{\veps,j}}$ then ${|\alpha_{n_k+j}-\alpha|<\veps}$, and consider the
largest integer ${k_\veps=\max_{j\in[1,\max j]}\{k_{\veps,j}\}}).$ Next
observe that
\vskip6pt\noi
{\narrower\narrower
if $\{n_k\}_{k\ge0}$ is of bounded increments, say with
$M=\sup_k({n_{k+1}-n_k})$, then
${\bigcup_{j\in[0,M]}\{n_k+j\}_{k\ge1}}=\{n\}_{n\ge1}.$
\vskip6pt}
\noi
So every positive integer $n$ is written as ${n=n_k+j}$ for some $n_k$ and
some ${j\in[1,M]}.$ Thus by the above displayed implication (with ${J=[1,M]}$
so that ${\max j=M}$), for every ${\veps>0}$ there exists an integer
$n_\veps={n_{k_\veps\!}+M}$
such that
$$
n\ge n_\veps\quad\limply\quad|\alpha_n-\alpha|<\veps.
$$
Hence ${\alpha_n\to\alpha}$ as ${n\to\infty}.$ The converse is trivial. \qed

\vskip6pt
As shown next, ${T^{n_k}x\wconv0}\limply{T^nx\wconv0}$ whenever
$\{T^{n_k}\}$ is \hbox{boundedly spaced}.

\vskip6pt\noi
{\bf Lemma 5.2.}
{\it $\kern-2.5pt$Suppose $\{T^{n_k}\}\kern-.5pt$ is a boundedly spaced
subsequence of $\{T^n\}.$ $\!$Take ${x\kern-.5pt\in\kern-.5pt\X}\kern-.5pt$.
$$
\hbox{\it If\/ $\{T^{n_k}x\}$ converges weakly to zero, then\/ $\{T^nx\}$
converges weakly to zero}\/.
$$

\proof
Take a vector ${x\in\X}$ and a subsequence $\{T^{n_k}\}_k$ of $\{T^n\}_n.$
Suppose $\{T^{n_k}x\}_k$ converges weakly to zero$.$ That is,
${\lim_kf(T^{n_k}x)}=0$ for every ${f\in\X^*}\!.$ Thus, in particular, for an
arbitrary $g$ in $\X^*$ (so that ${gT^j\!=T^{*j}}g$ lies in $\X^*$ for every
nonnega\-tive integer $j$), ${\lim_k(T^{*j}g)(T^{n_k}x)=0}.$ Now suppose
$\{T^{n_k}\}_k$ is a boundedly spaced subsequence of $\{T^n\}_n.\!$ Then
$\{(gT^j)(T^{n_k}x)\}_k$ is a boundedly spaced subsequence of the
$\FF$-valued sequence $\{(gT^j)(T^nx)\}_n$ for each ${j\ge0}.\!$ Since
${\lim_kg(T^{n_k+j}x)}={\lim_k(gT^j)(T^{n_k}x)}=0$ for every ${j\ge0}$
(i.e., for every ${gT^j\!=T^{*j}g\in\X^*}$), Proposition 5.1 ensures that
${\lim_ng(T^nx)}=0.$ As this holds for an arbitrary ${g\in\X^*}\!$, it follows
that $\lim_ng(T^nx)=0$ for every ${g\in\X^*}\!$, and so $\{T^nx\}_n$ converges
\hbox{weakly to zero}.                                                    \qed

\vskip6pt
Let ${T\kern-1pt\in\BX}$ be a weakly quasistable operator on a normed space
$\X.$ This means ${\liminf_n|f(T^nx)|=0}$ for every
${x\kern-.5pt\in\kern-.5pt\X}$, for every ${f\kern-.5pt\in\kern-.5pt\X^*}\!.$
Equivalently, for each ${x\in\X}$ there is at least one subsequence
$\{T^{n_k(x)}\}$ of $\{T^n\}$ for which ${\lim_k|f(T^{n_k(x)}x)|}=0$ for every
${f\in\X^*}\!$ (where, for each ${x\in\X}$, the sequence of integers
$\{n_k(x)\}$ depends on $x$ but not on $f).$ Such a subsequence
$\{T^{n_k(x)}\}$ of $\{T^n\}$ will be referred to as {\it a subsequence of
weak quasistability of\/ $T\kern-1pt$ for}\/ ${x\in\X}.$

\vskip6pt
Theorem 5.3 below is a central result.

\vskip6pt\noi
{\bf Theorem 5.3.}
{\it If\/ ${T\kern-1pt\in\BX}$ has a boundedly spaced subsequence of weak
quasistability for every\/ ${x\in\X}$, then\/ $T$ is weakly stable}\/.

\proof
Take an arbitrary ${x\in\X}.$ Suppose $T$ has a boundedly spaced subsequence
of weak quasistability $\{T^{n_k(x)}\}$ for such an $x.$ That is,
$\{T^{n_k(x)}x\}$ converges weakly to zero, and the subsequence
$\{T^{n_k(x)}\}$ of $\{T^n\}$ is boundedly spaced$.$ Then $\{T^nx\}$ converges
weakly to zero according to Lemma 5.2$.$ If this holds for every ${x\in\X}$,
then $T$ is weakly stable.                                                \qed

\vskip6pt
Proposition 4.3 and Theorem 5.3 together ensure the following nonimplication.
Take any ${x\in\X}.$ If ${\lim_kf(T^{n_k}x)=0}$ for every ${f\in\X^*}\!$, for
some subsequence $\{T^{n_k}\}$ of $\{T^n\}$, then it does not follow that
there exists a subsequence of bounded increments $\{m_k\}$ of $\{n\}$ for
which ${\lim_kf(T^{m_k}x)=0}$ for every ${f\in\X^*}$.

\section{\bf Application to weak l-sequential supercyclicity}

The orbit $\Oe_T(y)$ of a vector ${y\in\X}$ under an operator
${T\kern-1pt\in\kern-1pt\BX}$ is the set
$$
\Oe_T(y)={\bigcup}_{n\ge0}\{T^ny\}=\big\{T^ny\in\X\!:
\hbox{for every integer}\;n\ge0\big\}.
$$
Let $[x]=\span\{x\}$ denote the one-dimensional subspace of $\X$ spanned by a
sin\-gle\-ton $\{x\}$ at a vector ${x\in\X}.$ The projective orbit of a
nonzero vector ${y\in\X}$ under an operator ${T\in\BX}$ is the orbit
$\Oe_T([y])=\bigcup_{n\ge0}T^n([y])$ of the span of $\{y\}$:
$$
\Oe_T([y])=\big\{\alpha T^ny\in\X\!:\;
\hbox{for every}\;\alpha\in\FF\;\hbox{and every integer}\;n\ge0\big\}.
$$
A vector $y$ in $\X$ is a ({\it strongly}\/) {\it supercyclic vector}\/ for
$T$ if its projective orbit $\Oe_T([y])$ is dense in $\X.$ Since the norm
topology is metrizable, a nonzero vector $y$ in $\X$ is supercyclic for an
operator $T$ if and only if for every $x$ in $\X$ there exists an $\FF$-valued
sequence $\{\alpha_k\}_{k\ge0}$ such that for some sequence
$\{T^{n_k}\}_{k\ge0}$ with entries from $\{T^n\}_{n\ge0}$, the $\X$-valued
sequence $\{\alpha_kT^{n_k}y\}_{k\ge0}$ converges to $x$ (in the
\hbox{norm topology}):
$$
\alpha_kT^{n_k}y\sconv x
\;\;\hbox{for every}\;\;
x\in\X
\qquad(\hbox{i.e.,}\quad
\|\alpha_kT^{n_k}y-x\|\to0
\;\;\hbox{for every}\;\;
x\in\X).
$$
An operator $T$ is a ({\it strongly}\/) {\it supercyclic operator}\/ if it has
a supercyclic vector.

\vskip6pt
The weak version of the above convergence criterion leads to the notion of
weak l-sequential supercyclicity$.$ A nonzero vector $y$ in $\X$ is a
{\it weakly l-sequentially supercy\-clic vector}\/ for $T$ if for every $x$ in
$\X$ there is an $\FF$-valued sequence $\{\alpha_k\}_{k\ge0}$ such that, for
some sequence $\{T^{n_k}\}_{k\ge0}$ with entries from $\{T^n\}_{n\ge0}$, the
$\X$-valued sequence $\{\alpha_kT^{n_k}y\}_{k\ge0}$ converges weakly to $x.$
In other words, if
$$
\alpha_kT^{n_k}y\wconv x
\;\;\hbox{for every}\;\;
x\in\X
$$
\vskip-2pt\noi
$$
(\hbox{i.e.,}\quad
f(\alpha_kT^{n_k}y-x)\to0
\;\;\hbox{for every}\;\;
x\in\X,
\;\;\hbox{for every}\;\;
f\in\X^*).
$$
\vskip2pt\noi
Fix a weakly l-sequentially super\-cyclic vector $y.$ The $\FF$-valued
sequence $\{\alpha_k\}_{k\ge0}$ depends on $x$ (but not on $f$), and each
${\alpha_k=\alpha_k(x)}$ is nonzero whenever ${x\ne0}.$ Similarly, the
sequence $\{T^{n_k}\}_{k\ge0}$ also depends on $x$ (but not on $f$) --- and
so (as before) we sometimes write $T^{n_k}=T^{n_k(x)}.$ An operator $T$ is a
{\it weakly l-sequentially supercyclic operator}\/ if it has a weakly
l-sequentially supercyclic vector$.$

\vskip6pt\noi
{\bf Remark 6.1.}
(a)
Note that, as remarked in \cite{Kub2}, weak l-sequential supercyclicity is
not a topological notion$.$ In fact, as we have defined above, weak
l-sequential supercyclicity comes as a weak version of a norm convergence
criterion, but it does not mean denseness of the projective orbit in any
topology$.$ In particular, weak l-se\-quential supercyclicity does not mean
denseness of the projective orbit in weak topology or in weak sequential
topology, so weak l-sequential supercyclicity is not a weak topology notion$.$
Therefore, it is not the case to argue in terms of weak topology techniques
when dealing with weak l-sequential supercyclicity --- see \cite{Shk} for a
discussion along this line$.$ Also note that, according to the definition
of weak convergence in Section 2, weak l-sequential supercyclicity does not
mean convergence of the $\BX$-valued sequence
$\kern-1pt\{\alpha_kT^{n_k}\}\kern-1pt$ in the weak \hbox{operator topology
on $\kern-1pt\BX$}.

\vskip6pt\noi
(b)
Take a nonzero ${x\in\X}.$ If ${x\in\Oe_T([y])}$ for some vector $y$, then
$\{T^{n_k}\}_{k\ge0}$ can be viewed as a constant infinite sequence;
equivalently, as a single-entry finite subsequence of $\{T^n\}_{n\ge0}$ ---
i.e., $x=\alpha_0 T^{n_0}y$ for some ${\alpha_0\ne0}$ and some ${n_0\ge0}.$
In this case, the notions of ${\alpha_kT^{n_k}y\sconv x}$ and
${\alpha_kT^{n_k}y\wconv x}$, related to plain supercyclicity and weak
l-sequential supercyclicity, coincide$.$ In this case, convergence of the
constant or finite sequence $\{\alpha_kT^{n_k}y\}$ means that a limit is
eventually reached$.$ Therefore, in general, we use the expression
$$
\hbox{\it $\{T^{n_k}\}_{k\ge0}$ is a sequence with entries from\/
$\{T^n\}_{n\ge0}$}.
$$
However, if ${x\not\in\Oe_T([y])}$, this means
{\it $\{T^{n_k}\}_{k\ge0}$ is a subsequence of\/ $\{T^n\}_{n\ge0}$}.

\vskip6pt
For further weak forms of cyclicity see, e.g., \cite{BM}, \cite{CS},
\cite{DT}, \cite{Dug}, \cite{MS}, \cite{San1}, \cite{San2}$.$ Weak
l-sequential supercyclicity is the central theme in this section, and it has
been considered, for instance, in \cite{BCS}, \cite{BM1}, \cite{Kub1}, and
discussed in \cite{Shk}, \cite{KD1}, \cite{KD3}$.$ Observe that weak
l-sequential supercyclicity is not weak sequential supercyclicity$.$ For a
comparison between weakly l-sequentially supercyclicity and other notions of
weak cyclicity, including weak hypercyclicity, weak sequential supercyclicity,
and weak supercyclicity (plain weak cyclicity coincides with plain cyclicity)
see, for instance, \cite[pp.38,39]{Shk}, \cite[pp.259,260]{BM2},
\cite[pp.159,232]{GP}, \cite[pp.50,51,54]{KD1},
\hbox{\cite[pp.372,373,374]{KD2}.}

\vskip6pt
It has been known for a long time that {\it a power bounded supercyclic
operator on a normed space is strongly stable}\/ \cite[Theorem 2.2]{AB}$.$
Along this line, the following question has been posed in \cite{KD1}:
$$
\hbox{\it Is a power bounded weakly l-sequentially supercyclic operator
weakly stable $?$}
$$
This is a nontrivial question that, as far as we are aware, remains
unanswered$.$ The next results give partial answers to it.

\vskip6pt
Let ${Y_T\sse\X}$ denote the set of all weakly l-sequentially supercyclic
vectors $y$ for an operator $T$ on a normed space $\X$,
$$
Y_T=\big\{y\in\X\!:\;\hbox{for every}\;x\in\X\;\hbox{there exists}\;
\{\alpha_k\}\;\hbox{such that}\;x=w\hbox{\,-}{\lim}_na_kT^{n_k}y\big\},
$$
so that $T$ is weakly l-sequentially supercyclic if and only if
${Y_T\ne\void}.$ In this case, take an arbitrary ${y\in Y_T}.$ So, for every
${x\in\X}$, there is at least one sequence $\{T^{n_k}\}$ with entries from
$\{T^n\}$, and an associated scalar sequence $\{\alpha_k\}$, for which
${\lim_kf(\alpha_kT^{n_k}y)}\kern-.5pt=\kern-.51ptf(x)$ for every
${f\kern-.5pt\in\kern-.5pt\X^*}\!$. By Remark 6.1(b), if
${x\kern-.5pt\in\kern-.5pt\Oe_T([y])}$, \hbox{then the} above sequences can
be viewed as constant sequences, say ${x=\alpha_kT^{n_k}y}$ with
${\alpha_k\kern-1pt=\alpha_0}$ and ${n_k\kern-1pt=n_0}$ for every ${k\ge0}$,
where the constant sequence $\{n_k\}$ is not a subsequence of the nonnegative
integers, and so it makes no sense to talk of boundedly spaced subsequences in
this case$.$ We will, however, be concerned with boundedly spaced subsequences
$\{T^{n_k}\}$ of $\{T^n\}$ when $x$ lies in ${\X\backslash\Oe_T([y])}$ for an
arbitrary $y$ in $Y_T$.

\vskip6pt
Let ${T\!\in\kern-.5pt\BX}$ be a weakly l-sequentially supercyclic operator on
a normed \hbox{space $\kern-1pt\X\!.$} Let $y$ be any vector in $Y_T.$ For
each $x$ in ${\X\backslash\Oe_T([y])}$ there is at least one subsequence
$\{T^{n_k(x)}\}$ of $\{T^n\}$ for which $\lim_k\alpha_k(x)f(T^{n_k(x)}y)=f(x)$
for every ${f\in\X^*}\!$, where the $\FF$-valued sequence $\{\alpha_k(x)\}$
and the sequence of integers $\{n_k(x)\}$ depend on $x$ but not on $f.$ This
subsequence $\{T^{n_k(x)}\}$ of $\{T^n\}$ is referred to as {\it a subsequence
of weak l-sequential supercyclicity of\/ $T\kern-1pt$ for}\/
${x\in\X\backslash\Oe([y])}$ \hbox{with respect to\/ ${y\in Y_T}\!$.}}

\vskip6pt\noi
{\bf Theorem 6.2$.$}
{\it $\kern-2pt$Let $\kern-.5ptT\kern-2pt$ be a power bounded weakly
l-sequentially supercyclic \hbox{operator $\kern-1pt$on} a normed space
$\X\!.$ If there exists a boundedly spaced subsequence of weak l-sequential
supercyclicity for each\/ ${x\in\X\backslash\Oe_T([y])}$, for every\/
${y\in Y_T}$, then\/ $T$ is weakly stable}\/.

\proof
It has recently been proved in \cite[Corollary 4.3]{KD4} that
\vskip0pt\noi
\vbox{
$$
\hbox
{\it every power bounded weakly l-sequentially supercyclic}
$$
\vskip-8.5pt\noi
$$
                                                                 \eqno{(*)}
$$
\vskip-8pt\noi
$$
\hbox
{\it operator on a normed space is weakly quasistable}\/.
$$
}
\vskip0pt\noi
We sketch the proof below, both for sake of completeness and also because
part of its argument will be required later in this proof$.$

\vskip5pt\noi
Suppose $T$ a weakly l-sequentially super\-cyclic operator on a normed space
$\X$, and take an arbitrary ${y\in Y_T}.$ It was shown in
\cite[Theorem 4.1]{KD4} that if this $T$ is power bounded, then there is no
${f\in\X^*}$ such that ${\liminf_n|f(T^ny)|>0}.$ Consequently,
$$
{\liminf}_n|f(T^ny)|=0
\quad
\hbox{for every}\quad
y\in Y_T,
$$
for every ${f\in\X^*\!}.$ Also, we know that if ${Y_T\ne\void}$, then $Y_T$ is
dense in $\X$ \cite[Theorem 5.1]{Kub2}$.$ Then for an arbitrary ${x\in\X}$,
take (by density) any $Y_T$-valued sequence$\{y_k\}$ converging (in the norm
topology) to $x$, so that (since $T$ is power bounded)
$$
|f(T^nx)|\le|f(T^n(y_k-x))|+|f(T^ny_k)|
\le\|f\|\kern1pt({\sup}_m\|T^m\|)\kern1pt\|y_k-x\|+|f(T^ny_k)|
$$
for each integer ${n\ge0}$, for every ${f\in\X^*\!}.$ Hence, as
${\liminf_n|f(T^ny)|=0}$ for every $y$ in $Y_T$, ${\liminf_n|f(T^ny_k)|=0}$
for every $k$, and so (as $x$ was arbitrarily taken \hbox{from $\X$})
\vskip5pt\noi
$$
{\liminf}_n|f(T^nx)|=0
\quad
\hbox{for every}\quad
x\in\X,
$$
\vskip2pt\noi
for every ${f\in\X^*\!}$, proving the result in $(*).$

\vskip6pt\noi
Now we split the proof into two parts$.\kern-1pt$ Take ${y\in Y_T}$ and
${x\in\X\backslash\Oe_T([y])}$ arbitrarily$.$ Thus $T$ has a subsequence of
weak l-sequential supercyclicity $\{T^{n_k(x)}\}$ for $x$ (with respect to
$y$) so that, for some scalar sequence $\{\alpha_k(x)\}$,
\vskip5pt\noi
$$
{\lim}_k\alpha_k(x)f(T^{n_k(x)}y)=f(x)
\quad\hbox{for every}\quad
f\in\X^*\!,
$$
\vskip2pt\noi
where the above limit holds for every subsequence of
$\{\alpha_k(x)f(T^{n_k(x)}y)\}.$

\vskip5pt\noi
{\sc Part 1}.
Suppose $x$ is such that $\{\alpha_k(x)\}$ is unbounded$.$ So
$\limsup_k|\alpha_k(x)|=\infty.$ Since ${|f(x)|\in\RR}$, the above displayed
identity ensures that
\vskip5pt\noi
$$
{\liminf}_k|f(T^{n_k(x)}y)|=0
\quad\hbox{for every}\quad
f\in\X^*.
$$
\vskip2pt\noi
If $T$ has a boundedly spaced subsequence of weak l-sequential supercyclicity
for such an ${x\in\X\backslash\Oe_T([y])}$ and such a $y$, then take it so
that this same $\{T^{n_k(x)}\}$ is a boundedly spaced subsequence of weak
quasistability for $y$ as well$.$ So ${T^ny\wconv0}$ according to Lemma 5.2$.$
If this holds for every ${y\in Y_T}$, then ${T^ny\wconv0}$ for every
${y\in Y_T}.$ Since $Y_T$ is dense in $\X$, we get (as we saw above)
${T^nx\wconv0}$ for \hbox{every ${x\in\X}$.}

\vskip5pt\noi
{\sc Part 2}.
Suppose $x$ is such that $\{\alpha_k(x)\}$ is bounded$.$ If $T$ has a
boundedly spaced subsequence $\{T^{n_k(x)}\}$ of weak l-sequential
supercyclicity for such an ${x\in\X\backslash\Oe_T([y])}$ and such a $y$, then
take it so that the corresponding $\FF$-valued sequence $\{\alpha_k(x)\}$ is
boundedly spaced as well$.$ Since $\{\alpha_k(x)\}$ is bounded, there is a
convergent subsequence $\{\alpha_{k_j}\kern-1pt(x)\}$ of $\{\alpha_k(x)\}$ (by
the Bolzano--Weierstrass Theorem), which can be taken to be boundedly spaced
as well$.$ Hence $\{T^{n_{k_j}(x)}\}$ is a boundedly spaced subsequence of
weak l-sequential supercyclicity for $x$, so that
${\lim_j\alpha_{k_j}(x)f(T^{n_{k_j}(x)}y)}=f(x)$ for every ${f\in\X^*}\!$,
with ${\lim_j\alpha_{k_j}(x)=\alpha(x)}.$ If $T$ is power bounded, then it was
shown in \cite[Theorem 4.1(c)]{KD4} that
${\alpha(x)\lim_jf(T^{n_{k_j}(x)}y)}=f(x)$ for every ${f\in\X^*}\!.$ But this
implies that ${\liminf_j|f(T^{n_{k_j}(x)}y)|}=0$ (cf$.$
\cite[Theorem 4.1(d)]{KD4})$.$ Thus, as in Part 1, $\{T^{n_{k_j}}\}$ is a
boundedly spaced subsequence of weak quasistability for $y.$ Therefore the
same argument as in Part 1 ensures ${T^nx\wconv0}$ for every ${x\in\X}$.  \qed

\vskip6pt
Theorem 6.2 does not ensure the existence of a boundedly spaced subsequence of
weak l-sequential supercyclicity$.$ But as we have seen in the proof of
\hbox{Theorem 6.2} (cf$.$ \cite[Corollary 4.3]{KD4}), weak quasistability
holds under power boundedness$:$ {\it if\/ $T\kern-1pt$~is a weakly
l-sequentially supercyclic operator on a normed space\/ $\X$, then}
$$
\hbox{$T$ {\it power bounded}\/ $\;\limply\;\;T$ weakly quasistable}.
                                                             \leqno{\rm(a)}
$$
We show next that
$$ 
\hbox{$\{\alpha_k(x)\}$ {\it bounded for some}\/ $x\in\X\backslash\Oe_T([y])
\;\;\limply\;T$ weakly unstable}.                            \leqno{\rm(b)}
$$

\vskip6pt\noi
{\bf Theorem 6.3.}
{\it Let\/ $T\kern-1pt$ be a weakly l-sequentially supercyclic operator on a
normed space\/ $\X.\!$ Take an arbitrary vector\/ ${y\in Y_T}$ so that for
every vector\/ ${x\in\X\backslash\Oe_T([y])}$ there are a sequence of
scalars\/ $\{\alpha_k(x)\}$ and a subsequence\/ $\{T^{n_k(x)}\}\kern-1pt$ of\/
$\{T^n\}\kern-1pt$ such that\/ ${\alpha_k(x)f(T^{n_k(x)}y)}\to f(x)$ for
every\/ ${f\in\X^*}$ $\,($i.e., such that ${\alpha_k(T^{n_k}y)\wconv x}).$
\vskip6pt\noi
{\narrower\narrower
{\it If\/ for some\/ ${y\in Y_T}$ the sequence of scalars $\{\alpha_k(x)\}$
is bounded for some\/ ${x\in\X\backslash\Oe_T([y])}$, then\/ $T$ is not
weakly stable}\/.
\vskip0pt}
\vskip6pt\noi
{\it In other words, if\/ $T\kern-1pt$ is weakly stable, then all scalar
sequences\/ $\{\alpha_k\}$ are unbounded}\/.

\proof
Take an arbitrary ${y\in Y_T}$ and an arbitrary nonzero
${x_0\in\X\backslash\Oe_T([y])}.$ Although there is no injective linear
functional on a linear space of dimension greater than one (we are assuming
${\dim\X\kern-1pt>1})$, note that there exists ${f_0\in\X^*}$ such that
${f_0(x_0)\ne0}$ (since there is no ${0\ne z\in\X}$ for which ${f(z)=0}$
for every ${f\in\X^*}$ by the Hahn--Banach Theorem.) Thus, since $T$ is
weakly l-sequentially supercyclic, there exist a sequence of
scalars $\{\alpha_k(x_0)\}$ and a subsequence $\{T^{n_k(x_0)}\}\kern-1pt$ of
$\{T^n\}\kern-1pt$ such that
$$
|\alpha_k(x_0)|\kern1pt|f_0(T^{n_k(x_0)}y)|\to|f_0(x_0)|\in\RR\backslash\0.
$$
If $T$ is weakly stable, then ${\lim_n|f(T^ny)|=0}$, so that
${\lim_k|f(T^{n_k(x_0)}y)|=0}$, for \hbox{every} ${f\in\X^*}$; in particular,
for such an ${f_0\in\X^*}\!.$ Therefore the above displayed expression ensures
that the sequence $\{\alpha_k(x_0)\}$ is unbounded for such an arbitrary $x_0$
in ${\X\backslash\Oe_T([y])}.$ Consequently, if for some ${y\in Y_T}$ the
sequence of scalars $\{\alpha_k(x)\}$ is bounded for some
${x\in\X\backslash\Oe_T([y])}$, then the operator $T$ is not weakly stable.
                                                                          \qed

\vskip6pt
Take an arbitrary ${y\in Y_T}.$ Recall that an operator $T$ is said to have a
{\it boundedly spaced subsequence of weak l-sequential supercyclicity}\/ for
a given ${x\in\X\backslash\Oe([y])}$ if ${\alpha_k(x)f((T^{n_k(x)}y)\to f(x)}$
for every ${f\in\X^*}$ and the sequence $\{T^{n_k(x)}\}$ is boundedly spaced
(equivalently, the sequence of integers $\{n_k(x)\}$ is of bounded
increments).

\vskip6pt\noi
{\bf Corollary 6.4.}
{\it Let\/ $T$ be a power bounded weakly l-sequentially supercyclic operator
on a normed space\/ $\X.$ Then either
\begin{description}
\item{$\kern-4pt$\rm(a)}
$T$ does not have a boundedly spaced subsequence of weak l-sequential
supercyclicity for some ${x\in\X\backslash\Oe_T([y])}$ and some ${y\in Y_T}$,
\end{description}
\vskip-4pt\noi
or
\begin{description}
\item{$\kern-4pt$\rm(b)}
the sequence of scalars $\{\alpha_k(x)\}$ is unbounded for every
${x\in\X\backslash\Oe_T([y])}$, for every ${y\in Y_T}$.
\end{description}

\proof
If $T$ is weakly unstable, then Theorem 6.2 ensures the alternative (a), and
if $T$ is weakly stable, then Theorem 6.3 ensures the alternative (b),
exclusively.                                                             \qed

\vskip6pt
We close the paper with a question.
$$
\hbox{\it Is the Foguel operator in Proposition 4.3 weakly l-sequentially
supercyclic$\kern1pt$}?
$$


\vskip-10pt\noi
\bibliographystyle{amsplain}

\end{document}